\documentclass{conm-p-l}

\newtheorem{thm}{Theorem}
\newtheorem{lm}{Lemma}
\newtheorem{conj}{Conjecture}
\newcommand{\zset}{\mathbf{Z}} 
\newcommand{\rset}{\mathbf{R}} 
\newcommand{\bs}{{\bar s}}
\newcommand{\bbs}{{\bar{\bar s}}}

\numberwithin{equation}{section}

\begin{document}

\title[Blow-up in the Hamilton--Jacobi equation]%
{A Blow-Up Phenomenon in the Hamilton--Jacobi Equation\\ 
in an Unbounded Domain}

\author{K. Khanin}
\address{Isaac Newton Institute for Mathematical Sciences\\ 
20 Clarkson Rd\\ CB3~0EH Cambridge, UK}
\address{Heriot-Watt University\\ Edinburgh, UK}
\address{Landau Institute\\ Moscow, Russia}
\email{kk262@newton.cam.ac.uk} 

\author[D. Khmelev]{\protect\framebox{D. Khmelev}}
\address{Department of Mathematics\\
University of Toronto\\ M5S~3G3 Toronto, Ontario, Canada} 
\email{dkhmelev@math.toronto.edu}

\dedicatory{\parbox{0.7\textwidth}{Dmitry Khmelev died unexpectedly on
    24~October 2004.  Dmitry was a~very bright young mathematician and
    his tragic death at the age of 28 is a~big loss for the whole of
    the mathematical physics community.}}

\author{A. Sobolevski{\u\i}}
\address{Department of Physics\\
M.V.~Lomonossov Moscow State University\\ 
Moscow, Russia}
\curraddr{Laboratoire G.D. Cassini\\
Observatoire de la C{\^o}te d'Azur\\
BP~4229, 06304 Nice, France}
\email{ansobol@obs-nice.fr} 
\thanks{The work of AS was supported by the INTAS fellowship YSF2001:172 and the RFBR grant 02--01--1062}
 

\begin{abstract}
  We construct an example of blow-up in a flow of min-plus linear
  operators arising as solution operators for a Hamilton--Jacobi
  equation $\partial S/\partial t+|\nabla S|^\alpha/\alpha+U(x,t)=0$,
  where $\alpha>1$ and the potential $U(x,t)$ is uniformly bounded
  together with its gradient.  The construction is based on the fact
  that, for a suitable potential defined on a time interval of
  length~$T$, the absolute value of velocity for a Lagrangian
  minimizer can be as large as $O\bigl((\log
  T)^{2-2/\alpha}\bigr)$. We also show that this growth estimate
  cannot be surpassed. Implications of this example for existence of global
  generalized solutions to randomly forced Hamilton--Jacobi or Burgers
  equations are discussed.
\end{abstract}

\maketitle

\section{Introduction}

In this paper we present an example of blow-up in a flow of min-plus
linear integral operators arising as solution operators for a class of
Hamilton--Jacobi equations. As we shall see, existence of such blow-up
has interesting consequences for the application of idempotent
functional analysis to stochastic partial differential equations.

\subsection{}
Consider the inviscid Burgers equation in the $d$-dimensional space
\begin{equation}
  \label{burgers}
  \frac{\partial u}{\partial t}+(u\cdot\nabla)u=-\nabla U(x,t),
\end{equation}
where $u(x,t)=(u_1(x,t),u_2(x,t),\dots,u_d(x,t))$ is a potential
velocity field, so that $u(x,t)=\nabla S(x,t)$. The potential $S(x,t)$
must satisfy the Hamilton--Jacobi equation
\begin{equation}
  \label{HJEsquare}
  \frac{\partial S}{\partial t}+\frac12\left|\nabla S\right|^2+U(x,t)=0.
\end{equation}
Here and below, $\nabla$ denotes the vector of derivatives with
respect to components of the vector $x\in\rset^d$.

It is well-known that the Cauchy problems for nonlinear equations
\eqref{burgers} and~\eqref{HJEsquare} fail to have global in time
classical solutions: they develop infinite velocity gradients in
finite time. There exist several ways to extend solutions beyond
formation of such singularities in a suitable generalized sense,
allowing for discontinuities of velocities
\cite{H50,L82,CL83,S95,KM97}. Under an additional stability hypothesis, 
all of them become essentially equivalent
(see, e.g., the paper \cite{R03} in the present volume), and the
corresponding solutions admit an explicit representation in terms of
the Lax--Ole{\u\i}nik variational principle.

Namely, a generalized solution to a Cauchy problem for the
Hamilton--Jacobi equation~\eqref{HJEsquare} with the initial condition
$S(x,0)=S_0(x)$ has the form
\begin{equation}
  \label{laxoleinik}
  S(x,t)=\inf_{\gamma(t)=x}
  \bigl(A_{0,t}[\gamma]+S_0\bigl(\gamma(0)\bigr)\bigr),
\end{equation}
where the action functional $A_{\cdot,\cdot}[\cdot]$ is given by
\begin{equation}
  \label{defaction}  
  A_{t_1,t_2}[\gamma]\equiv\int_{t_1}^{t_2} L(\dot\gamma(s),\gamma(s),s)\,ds
\end{equation}
for any $t_1$ and~$t_2$ with $t_1<t_2$, the Lagrangian has the form
$L(v,x,t)=|v|^2/2-U(x,t)$, and the infimum is taken over all
absolutely continuous trajectories~$\gamma(\cdot)$ defined over
$[0,t]$ and satisfying $\gamma(t)=x$. Define further
\begin{equation}
  \label{defkernel}
  A_{t_1,t_2}(y,x)=\inf_{\gamma(t_1)=y,\;\gamma(t_2)=x}
  A_{t_1,t_2}[\gamma].
\end{equation}
Under mild conditions on the Lagrangian, this infimum, as well as the
infimum in~\eqref{laxoleinik}, is attained at a
trajectory $\gamma^{y,x}_{t_1,t_2}\colon [t_1,t_2]\to\rset^d$ (see, e.g.,
\cite{F01}); below we call such trajectories Lagrangian minimizers. The
solution to the Cauchy problem for the Burgers
equation~\eqref{burgers} on the time interval $[0,t]$ with the initial
condition $u(x,0)=\nabla S_0(x)$ is then given by
$u(x,t)=\dot\gamma^x_{0,t}(t)$, where $\gamma^x_{0, t}$ is a
Lagrangian minimizer corresponding to the minimum in the right-hand
side of
\begin{equation}
  \label{laxoleinikkernel}
  S(x,t)=T_{0,t}S_0(x)\equiv\min_y\bigl(A_{0,t}(y,x)+S_0(y)\bigr).
\end{equation}

For the purposes of the present paper, the Lax--Ole{\u\i}nik formula
\eqref{laxoleinik} or~\eqref{laxoleinikkernel} constitutes a
sufficient replacement for definitions of generalized solutions. Note
that in its form~\eqref{laxoleinikkernel}, the Lax--Ole{\u\i}nik
formula becomes a min-plus integral operator representation of a
solution.  The solution operators $T_{\cdot,\cdot}$ form a flow, i.e.,
they satisfy $T_{t_2,t_3}T_{t_1,t_2}=T_{t_1,t_3}$ for any
$t_1<t_2<t_3$; however, this flow is not $t$-translation invariant unless
$U(x,t)$ does not depend on time.

We note that the duality between representations of solutions in terms
of the value function $S(x,t)$ or minimizers $\gamma_{0,t}^x$ is more
than a heuristic relation; when one relaxes the action minimization
problem in the spirit of Kantorovich, allowing measure-valued
solutions instead of classic minimizing curves, the function~$S(x,t)$
becomes the dual variable in a correponding infinite-dimensional
linear program (see, e.g., \cite{M89,EG02}).

\subsection{}
Our interest in solution operators of the
form~\eqref{laxoleinikkernel} is motivated by the theory of global
(time-stationary) viscosity solutions in the case of randomly forced
inviscid Burgers and Hamilton--Jacobi equations, which was developed
recently in \cite{EKMS00}, \cite{IK03} and~\cite{GIKP03}. The crucial
role in the construction of this global solution is played by
Lagrangian minimizers $\gamma^x_t$ defined over a semi-infinite
time interval $(-\infty,t]$: namely, a global solution to the random
forced inviscid Burgers is given by $u(x,t)=\dot\gamma^x_{t}(t)$. To
prove that such semi-infinite minimizers exist, one has to take a limit
as $T\to\infty$ for minimizers $\gamma^x_{t-T,t}$ defined on
finite time intervals of the form $[t-T,t]$. Existence of this limit
follows from a uniform bound on the absolute value of a velocity
$|\dot\gamma^x_{t-T,t}(t)|$, which thus becomes the central problem
for the theory.

Observe first that the velocity of a minimizer is uniformly bounded if
the state space of the Lagrangian system is a compact
manifold~$M$. Indeed, in this case the displacement of a minimizer for
any time interval is bounded by the diameter of the manifold, so
action minimizing trajectories cannot have large velocities. The
simplest example is given by the $d$-dimensional torus
$\rset^d/\zset^d$. Hence, the uniform bound on velocities holds in the
case of $\zset^d$-periodic potential $U(x,t)$, satisfying
$U(x+k,t)=U(x,t)$ for all $k\in\zset^d$. It turns out that, for the
randomly forced Burgers equation on a compact manifold, a unique
global solution~$u(x,t)$ exists with probability~1. In fact the whole
theory is developed at the moment only in the case of compact
manifolds, where the bound on velocities can be easily proved. At
present almost nothing is known about global solutions in the case of
$\rset^d$ (however see \cite{HK03} for some results and discussions).

In the case of non-periodic potentials one can imagine a situation
where a minimizer spends almost all its time in a very favourable part
of $\rset^d$ which may lie far away from its prescribed endpoint~$x$,
and then goes very quickly to~$x$. Such scenario will lead to a large
terminal velocity at point~$x$ which might depend on the time
interval where minimization is performed. There are two cases, however, when such behaviour is
impossible. The first one corresponds to the autonomous bounded potential:
$U(x,t)=U(x)$, for which the energy
\begin{equation}
  \label{energy}  
  H(p,x,t)=\max_{v\in\rset^d}\bigl(p\cdot v-L(v,x,t)\bigr)
  =\frac{|p|^2}{2}+U(x,t)
\end{equation}
is conserved and the velocity of any Lagrangian trajectory is
uniformly bounded if this trajectory is at rest at the initial moment
of time. Since all minimizers are Lagrangian trajectories, the bound
on their velocities follows immediately.

The second case corresponds to a potential $U(x,t)$ that depends on
time periodically. Here the situation is more delicate. It is not true
anymore that the velocities of Lagrangian trajectories are
bounded. Moreover, it was shown recently by J.~Mather that Lagrangian
trajectories can be accelerated by a periodic potential to an
arbitrary large velocity even on a compact manifold. However, A.~Fathi was able to show with
methods developed in~\cite{F01} that the velocities of \emph{minimizing} trajectories are
still bounded; his elegant unpublished proof is recalled in
Appendix~\ref{a:fathi} below.

The examples constructed in this paper show that for special
potentials $U(x,t)$ the velocity of a minimizer may be arbitrarily
large; in fact, one can construct a potential $U(x,t)$ defined for all
$t<0$ that accelerates minimizers to infinite velocitites. Because of
this blow-up in velocity, for such potentials even generalized global
solutions do not exist. The simple remarks we just made demonstrate
that for this blow-up effect it is crucial that the system be defined
on an unbounded manifold (say $\rset^d$) and the potential $U(x,t)$
depend on time non-periodically. Implications of our examples to the
existence of global solutions in the randomly forced case is discussed
in the conclusion to this paper.

\subsection{}
We pass now to precise formulation of our results. Below we consider
not~\eqref{HJEsquare} but a more general Hamilton--Jacobi equation
\begin{equation}
  \label{HJE}
  \frac{\partial S}{\partial t}+H(\nabla S,x,t)=0,
\end{equation}
where the Hamiltonian has the form
\begin{equation}
  \label{defhamiltonian}
  H(p,x,t)=\frac1\alpha|p|^\alpha+U(x,t).
\end{equation}
The corresponding Lagrangian system has the Lagrangian
\begin{equation}
  \label{Lag}  
  L(v,x,t)= \frac{1}{\beta}|v|^\beta- U(x,t),
\end{equation}
where $\alpha^{-1}+\beta^{-1}=1$. Suppose that $\alpha,\beta>1$ and the
potential $U(\cdot,t)$ is a $C^1$~function of~$x$ for any~$t$,
uniformly bounded together with its spatial derivative:
\begin{equation}
  \label{Bound}  
  0\leq U(x,t)\leq C,\;
  |\nabla U(x,t)|\leq C,\quad 
  x\in\rset^d, t\in\rset.
\end{equation}

Let the trajectory $\gamma^x_{t_1,t_2}\colon [t_1,t_2]\to\rset^d$ be a
(not necessarily unique) Lagrangian minimizer for the action~$A_{t_1,t_2}$,
satisfying the conditions $\dot\gamma^x_{t_1,t_2}(t_1)=0$,
$\gamma^x_{t_1,t_2}(t_2)=x$. Note for future references that, under
the above conditions on Lagrangian, $\gamma^x_{t_1,t_2}$ is a
classical solution of the Euler--Lagrange equation
\begin{equation}
  \label{EL}
  \frac{d}{dt}(\dot\gamma(t)|\dot\gamma(t)|^{\beta-2}) = -\nabla U
\end{equation}
(see, e.g., \cite{F01}), where the dot notation stands for the
ordinary derivative with respect to time variable.

\begin{thm}
  \label{upperbound} 
  There exists $K_1=K_1(C,\beta)>0$ such that for any $[t_1,t_2]$ with
  large enough $T\equiv t_2-t_1$ and any $x\in\rset^d$
  \begin{equation}
    \label{Upperbound}
    |\dot\gamma^x_{t_1,t_2}(t_2)|\leq K_1(\log T)^{2/\beta}.
  \end{equation}
\end{thm}

\begin{thm}
  \label{lowerbound}
  There exists $K_2=K_2(C,\beta)>0$ such that for any $[t_1,t_2]$ with
  large enough $T\equiv t_2-t_1$ and any $y\in\rset^d$ there is a
  potential $U(\cdot,t)$, defined on the time interval $[t_1,t_2]$ and
  satisfying~\eqref{Bound}, such that
  \begin{equation}
    |\dot\gamma^x_{t_1,t_2}(t_2)|\geq\frac{K_2(\log T)^{2/\beta}}{2^{\beta/(\beta-1)}}
  \end{equation}
  for any $x$ with $|x-y|\leq R_T\equiv\frac{K_2}{2}(\log T)^{2/\beta}$.
\end{thm}

Later on constants in Theorems \ref{upperbound} and~\ref{lowerbound}  will be given explicit
expression in terms of the parameters $C$ and~$\beta$.

\begin{thm}
  \label{semiinf}
  There exists a potential $U(x,t)$, defined for all $t<0$ and
  satisfying~\eqref{Bound}, such that for all $x \in \rset^d$
  \begin{equation}
    \limsup_{t\to-\infty}|\dot\gamma^x_{t,0}(0)|=\infty.
  \end{equation}
  Moreover, the potential $U(x,t)$ may be chosen continuous in time.
\end{thm}

The paper is organized as follows. Theorem~\ref{upperbound} is proved
in Section~\ref{s:upperbound}. Theorems \ref{lowerbound}
and~\ref{semiinf} are proved in Section~\ref{s:lowerbound}. In
Section~\ref{s:conclusion}, we make concluding remarks and indicate
several directions in which one can generalize the results of the
present paper. In Appendix~\ref{a:g_T} we give the technical proof of
Lemma~\ref{t:g_T}, deferred from the main
text. Appendix~\ref{a:fathi}, included for completeness, contains
A.~Fathi's argument that rules out blow-up if the potential $U(x,t)$
is periodic in time.

To simplify notation we denote below the minimizer
$\gamma^x_{t_1,t_2}$ by~$\gamma^x$ and assume that all constants may
have implicit dependence on the parameters $C$ and~$\beta$. For
convenience we introduce a positive variable $s=t_2-t$ for $t \in
[t_1,t_2]$ and denote by $w(s)$ the absolute value of the average
velocity over $[0,s]$:
\begin{equation}
  \label{defws}
  w(s)\equiv\frac{|\gamma^x(t_2)-\gamma^x(t_2-s)|}{s}.
\end{equation}

\section{Proof of the upper bound on velocity}
\label{s:upperbound}

Before giving the proof of Theorem~\ref{upperbound} in full
generality, we observe that it becomes particularly simple in the case
of $\beta=2$. Fix a time interval $[t_1,t_2]$ and a minimizer
$\gamma^x$ with final position $\gamma^x(t_2)=x$. Take $s_1$ and~$s_2$
with $0\leq s_1\leq s_2\leq T$, where $T\equiv t_2-t_1$, and suppose that
the absolute value of the average velocity of the minimizer increases
from $w_2\equiv w(s_2)$ to $w_1\equiv w(s_1)$ over the time interval
$[t_2-s_2,t_2-s_1]$.% (see Fig.~\ref{f:w(s_2)}).
%\begin{figure}
%  \centerline{\psfig{file=fig01.eps,width=6cm}}
%  \caption{Construction used in the proof of the upper bound on
%  velocity (shown is the one-dimensional case).}
%  \label{f:w(s_2)}
%\end{figure}

Observe that minimization of the action allows to control the increase
in the average velocity:
\begin{equation}
  \label{w(s_2)quad}
  1+\frac{(w_1-w_2)^2}{2C}\leq\frac{s_2}{s_1}.
\end{equation}
To see this, note that
\begin{equation}
  \label{w(s_2)t1}
  \begin{split}
    A_{t_2-s_2,t_2}[\gamma^x]&=
      A_{t_2-s_2,t_2-s_1}[\gamma^x]+A_{t_2-s_1,t_2}[\gamma^x]\\
    &\geq \frac12\bigl(s_1w_1^2+\frac1{s_2-s_1}(s_1w_1-s_2w_2)^2\bigr)-Cs_2,
  \end{split}
\end{equation}
where to estimate the action we use~\eqref{Bound} and Jensen's
inequality, taken in the form
\begin{equation}
  \label{estimate}  
  \int_{t'}^{t''} |\dot\gamma(t)|^\beta\, dt \geq 
  (t''-t')^{1-\beta}|\gamma(t'')-\gamma(t')|^\beta
\end{equation}
for $\beta>1$ and an arbitrary $C^1$ curve $\gamma(t)\colon
[t',t'']\to\rset^d$. On the other hand, consider a trajectory
$\gamma(t)$, $t\in[t_2-s_2,t_2]$, that has the same endpoints
as~$\gamma^x$ but keeps constant velocity, which is equal to~$w_2$. By
action minimization and \eqref{Bound},
\begin{equation}
  \label{w(s_2)t2}
  A_{t_2-s_2,t_2}[\gamma^x]\leq A_{t_2-s_2,t_2}[\gamma]
  \leq\frac12s_2w_2^2.
\end{equation}
Combining \eqref{w(s_2)t1} and~\eqref{w(s_2)t2}, after some simple
algebra we arrive at~\eqref{w(s_2)quad}.

The meaning of inequality~\eqref{w(s_2)quad} is that increasing the
absolute value of the average velocity in \emph{arithmetic}
progression requires a \emph{geometric} progression in time steps.
Therefore the largest possible increase over a time interval of
length~$T$ is proportional to $\log T$. The desired
bound~\eqref{Upperbound} on the terminal velocity~$\dot\gamma^x(t_2)$
may now be inferred from (i)~the observation that the smaller is the
time interval, the closer are the absolute values of average and
terminal velocity, and (ii)~the boundedness of the average velocity
$w(T)$ at the earliest time moment~$t_1=t_2-T$, which we prove in a
separate lemma for future reference.

\begin{lm}
  \label{t:w(T)}
  $w(T)\leq (C\beta)^{1/\beta}$.
\end{lm}

\begin{proof}
Using~\eqref{estimate} and~\eqref{Bound}, it is easy to see that
\begin{equation}
  A_{t_1,t_2}[\gamma^x]\geq(T/\beta)(w(T))^\beta-CT.
\end{equation}
On the other hand, the action of the curve $\gamma(t)=x$ for all
$t\in[t_1,t_2]$, satisfies the estimate $A_{t_1,t_2}[\gamma]\leq
0$. Since $A_{t_1,t_2}[\gamma^x]\leq A_{t_1,t_2}[\gamma]$, we have
$w(T)\leq (C\beta)^{1/\beta}$.
\end{proof}

Turning now to the proof of Theorem~\ref{upperbound}, we start
with two auxiliary results. The first lemma extends
inequality~\eqref{w(s_2)quad} to the case of general $\beta>1$. 

\begin{lm}
\label{t:w(s_2)}
For $0\leq s_1\leq s_2\leq T$ denote
\begin{equation}
  \label{defDelta}
  w_1\equiv w(s_1),\quad
  w_2\equiv w(s_2),\quad
  \Delta\equiv w_1-w_2=\xi w_1^{(2-\beta)/2}
\end{equation}
and assume $0<\Delta<w_1$. There exists $W=W(\xi)>0$ such that if $w_1>W$,
then
\begin{equation}
  \label{s_2/s_1}  
  1+\frac{\xi^2(\beta-1)}{3C}\leq\frac{s_2}{s_1}.
\end{equation}
\end{lm}

\begin{proof}
Using~\eqref{estimate} and~\eqref{Bound}, we get
\begin{equation}
  \label{bound1}
  \begin{array}{r@{\,}c@{\,}l}
    A_{t_2-s_2,t_2}[\gamma^x]&=&
      A_{t_2-s_1,t_2}[\gamma^x]+A_{t_2-s_2,t_2-s_1}[\gamma^x]\\[1ex]
    &\geq&\displaystyle\frac{s_1w_1^\beta}{\beta}+
     \frac{(s_2-s_1)^{1-\beta}}{\beta}|\gamma^x(t_2-s_1)-\gamma^x(t_2-s_2)|^\beta-Cs_2\\[2ex]
    &\geq&\displaystyle\frac 1\beta\bigl(s_1w_1^\beta+
     (s_2-s_1)^{1-\beta}|s_2w_2-s_1w_1|^\beta\bigr)-Cs_2.
  \end{array}
\end{equation}
Denote by $\gamma(t)$, $t\in[t_2-s_2,t_2]$, the trajectory of
a point which moves with constant velocity from
$(\gamma^x(t_2-s_2),t_2-s_2)$ to $(\gamma^x(t_2),t_2)$. Since
\begin{equation}
  \label{bound2}
  A_{t_2-s_2,t_2}[\gamma]\leq\frac{s_2w_2^\beta}{\beta}
  =\frac{s_2(w_1-\Delta)^\beta}{\beta}
\end{equation}
and $A_{t_2-s_2,t_2}[\gamma^x]\leq A_{t_2-s_2,t_2}[\gamma]$,
inequalities \eqref{bound1} and~\eqref{bound2} imply
\begin{equation}
  s_1w_1^\beta+(s_2-s_1)\left|w_1-\frac{s_2}{s_2-s_1}\Delta\right|^\beta
  -C\beta s_2\leq s_2(w_1-\Delta)^\beta.
\end{equation}
With the notation $\sigma\equiv s_2/(s_2-s_1)$, this inequality is
equivalent to
\begin{equation}
  \label{temp0}
  \bigl|1-(\sigma\Delta/w_1)\bigr|^\beta\leq
  1+\sigma\Bigl((1-\bigl(\Delta/w_1)\bigr)^\beta-1+C\beta w_1^{-\beta}\Bigr).
\end{equation}
Using in the right-hand side of this inequality the Taylor expansion
$(1-z)^\beta=1-\beta z+\frac{\beta(\beta-1)}{2}z^2 (1-\theta(z)
z)^{\beta-2}$ with $\theta(z)\in[0,1]$, we get:
\begin{equation}
 \label{bound3}
  \left|1-\frac{\sigma\Delta}{w_1}\right|^\beta\leq
  1-\frac{\beta\sigma\Delta}{w_1}\left(1
  -\frac{(\beta-1)}{2}\frac{\Delta}{w_1}
  \left(1-\theta(\Delta/w_1)\frac\Delta{w_1}\right)^{\beta-2}\!\!\!\!-
  \frac{C}{\Delta w_1^{\beta-1}}\right).
\end{equation}
Since $\Delta/w_1$ and $1/(\Delta w_1^{\beta-1})$ are both of class $O(w_1^{-\beta/2})$ for
fixed $\xi$, the value of the
largest parenthesis in the right-hand side of~\eqref{bound3} lies between
$2(1+\beta)^{-1}$ and~$1$ if $w_1>W$ with a suitably
large~$W=W(\xi)$. Since the left-hand side of~\eqref{bound3} is
nonnegative, this implies
\begin{equation}
  \label{bound4}
  \frac{\sigma\Delta}{w_1}\le\frac{\beta+1}{2\beta}<1  
\end{equation}
and enables us to
use the same expansion in the left-hand side of~\eqref{bound3}. After
some cancellations this leads to the inequalities
\begin{equation}
  \label{bound5}
  \begin{array}{r@{\;}c@{\;}l}
    \sigma&\leq&\displaystyle\frac1{(1-\theta(\sigma\Delta/w_1)\sigma\Delta/w_1)^{\beta-2}}
  \left((1-\theta(\Delta/w_1)\Delta/w_1)^{\beta-2}+\frac{2C}{\xi^2(\beta-1)}\right)\\[3ex]
  &\leq&\displaystyle\max\left\{1,\left(\frac{2\beta}{\beta-1}\right)^{\beta-2}\right\}\left(1+\frac{2C}{\xi^2(\beta-1)}\right),
  \end{array}
\end{equation}
where the last line follows
from~\eqref{bound4} if $w_1>W$. The second of these inequalities says that
for $\xi$ fixed, $\sigma$ is bounded above uniformly in~$w_1$. Using this upper estimate on~$\sigma$
and enlarging $W$ if necessary, we can ensure that for $w_1>W$ the
parentheses containing $\theta$ in the first of inequalities~\eqref{bound5} are
arbitrarily close to unity, and therefore
\begin{equation}
  \label{last}
  \sigma\leq 1+\frac{3C}{\xi^2(\beta-1)},
\end{equation}
which implies~\eqref{s_2/s_1}.
\end{proof}

Note that in~\eqref{last}, as well as in~\eqref{s_2/s_1}, the
constant~$3$ may be replaced by any number greater than~$2$.

Using inequality~\eqref{s_2/s_1}, one can replace the arithmetic
progression in the~$w$ variable, suggested by
bound~\eqref{w(s_2)quad}, by a more general sequence that still leads
to a power-law estimate in~$\log T$ for the average velocity.  The
following lemma, employed several times throughout this paper, shows
that such estimate allows to control the terminal velocity
$\dot\gamma^x(t_2)$.
\begin{lm}
  \label{t:w_1}
  If $w(s)\leq (2Cs)^{1/(\beta-1)}$ then $|\dot\gamma^x(t_2)|\leq
  (3Cs)^{1/(\beta-1)}$. If $w(s)>(2Cs)^{1/(\beta-1)}$ then
  \begin{equation}
    \label{averagedot}
    (1/2)^{1/(\beta-1)}w(s)\leq|\dot\gamma^x(t_2)|\leq(3/2)^{1/(\beta-1)}w(s).
  \end{equation}
\end{lm}

\begin{proof}
The minimizer $\gamma^x(t)$ satisfies the Euler--Lagrange
equation~\eqref{EL}. This together with~\eqref{Bound} implies
\begin{equation}
  \label{EL1}  
  \bigl||\dot\gamma^x(t')|^{\beta-1}-|\dot\gamma^x(t'')|^{\beta-1}\bigr|\leq Cs
\end{equation}
for all $t',t'' \in [t_2-s,t_2]$. Since the Lagrangian~\eqref{Lag} is
strictly convex, $\gamma^x(t)$ is a $C^1$ curve, and there exists $t^*
\in [t_2-s,t_2]$ such that $|\dot\gamma^x(t^*)|=w(s)$. It follows
from~\eqref{EL1} written for $t'=t_2$ and $t''=t^*$ that 
\begin{equation}
  w(s)^{\beta-1}-Cs\leq|\dot\gamma^x(t_2)|^{\beta-1}\leq
  w(s)^{\beta-1}+Cs,
\end{equation}
which implies the statement.
\end{proof}

\begin{proof}[Proof of Theorem~\ref{upperbound}]
Somewhat departing from notation of Lemma~\ref{t:w(s_2)}, denote
$w_1\equiv|\gamma(t_2-1)-\gamma(t_2)|$, $\Delta\equiv
w_1^{(2-\beta)/2}$, $\bar W\equiv\sup\{\, W(\xi) \mid
2^{(2-\beta)/2}\leq\xi\leq 1\,\}$. Suppose that
$w_1>\max\{2(C\beta)^{1/\beta}, 2\bar W,(2C)^{1/(\beta-1)}\}$; otherwise the statement is trivially satisfied for large enough~$T$. 

Denote $s_0\equiv 1$. Since
$w(T)\leq (C\beta)^{1/\beta}$ by Lemma~\ref{t:w(T)} and $w(s)$ is a
continuous function, we can choose an increasing sequence of time instants $s_0<s_1<\dots<s_n$ such that $w(s_0)=w_1$, $w(s_1)=w_1-\Delta$, $\ldots$, $w(s_n)=w_1-n\Delta$ and $n=[w_1/(2\Delta)]$, where $[\cdot]$
stands for the integer part.
Denote $\xi_i\equiv\Delta w(s_i)^{(\beta-2)/2}$. Since
$w_1/2 \leq w(s_i) \leq w_1$ for $0\leq i\leq n$, all $\xi_i$ satisfy
the inequalities $2^{(2-\beta)/2}\leq \xi_i\leq1$ and therefore, by the choice of~$w_1$, all
$w(s_i)$ satisfy the condition of Lemma~\ref{t:w(s_2)}: $w(s_i)>\bar
W\geq W(\xi_i)$. Hence
\begin{equation}
  T\geq s_n=\prod_{i=1}^n\frac{s_i}{s_{i-1}}\geq 
  \left(1+\frac{2^{2-\beta}(\beta-1)}{3C}\right)^n.
\end{equation}
It follows that $n\leq\widetilde K\log T$, where
$\widetilde
K{=}\bigl(\log(1+\frac{2^{2-\beta}(\beta-1)}{3C})\bigr)^{-1}$.
Taking into account that $n=[w_1^{\beta/2}/2]$, we get $w_1\leq\max\{2(C\beta)^{1/\beta},2\bar W,(2C)^{1/(\beta-1)},
(2\widetilde K\log T+2)^{2/\beta}\}$ for any $T>0$. The statement now 
follows from Lemma~\ref{t:w_1}.
\end{proof}

\section{Construction of accelerating potentials}
\label{s:lowerbound}

Recall that $[t_1,t_2]$ is a fixed time interval with $t_2-t_1=T$.
To prove Theorems \ref{lowerbound} and~\ref{semiinf}, it is enough to
construct in this time interval an example of a potential that depends
only on one spatial coordinate. Hence, without loss of generality, we
may assume $d=1$, $x \in \rset$.

Observe that setting $s_0$ equal to $s$ instead of~$1$ in the proof of
Theorem~\ref{upperbound} gives for the average velocity of a minimizer at time
$t_2-s$ the bound~$O\bigl((\log(T/s))^{2/\beta}\bigr)$, which can be turned into a similar bound on~$\dot\gamma(t_2-s)$ by an argument analogous to that of Lemma~\ref{t:w_1}. For $s\in[0,T]$ and any $K>0$,
define
\begin{equation}
  \label{defgT}
  g_T(s)\equiv K\int_0^s(\log(T/u))^{2/\beta}\,du.
\end{equation}
Intuitively, this formula means that the trajectory $-g_T(t_2-t)$ has the ``largest velocity
possible'' for a minimizer at all times $t\in[t_1,t_2]$, up to the constant factor~$K$; accelerating
potentials constructed below confine minimizers to lie as close to
this trajectory as possible.

Before starting the proofs of Theorems \ref{lowerbound}
and~\ref{semiinf}, we collect here some properties of the
function~$g_T(\cdot)$ for future references.
\begin{lm}
  \label{t:g_T}
  Let $0\leq s\leq T$. Then
  \begin{gather}
    \label{integr}
    \int_0^s\frac1\beta(\dot g_T(u))^\beta\,du=
    \frac{K^\beta}{\beta}s\bigl((\log(T/s))^2+2\log(T/s)+2\bigr),\\
    \label{residue}
    g_T(s)=Ks\bigl(\log(T/s))^{2/\beta}
      \left(1+\frac2{\beta\log(T/s)}
        +\frac{2(2-\beta)}{\beta^2}r\bigl(\log(T/s)\bigr)\right),
  \end{gather}
  where $0\leq r(z)\leq z^{-2}$ for~$z>0$, and
  \begin{equation}
    \label{Main}
    \int_0^s\frac1\beta(\dot g_T(u))^\beta\,du-
      \frac{s^{1-\beta}}{\beta}\bigl(g_T(s)\bigr)^\beta<
      \frac{4K^\beta s}{\beta}.
  \end{equation}
  If $T>T_0$ for a suitable $T_0$ and $3<s\leq T$,
  then there exists $\bar M>0$ such that
  \begin{equation}
    \label{s_2}
    \frac{(g_T(s)-g_T(1))^\beta}{(s-1)^{\beta-1}}-
      \frac{(g_T(s))^\beta}{s^{\beta-1}}\leq 
      \bar M K^\beta(\log T)^2.
  \end{equation}
\end{lm}
The proof is postponed to Appendix~\ref{a:g_T}.

\subsection{Proof of Theorem~\ref{lowerbound}}
For any $y\in\rset$ define on the time interval $[t_1,t_2]$ a
potential
\begin{equation}
  \label{defpot1}
  U(x,t)\equiv U_C(x-y+g_T(t_2-t)),
\end{equation}
where $U_C(\cdot)$ is a $C^1$ function that satisfies the conditions
$0\leq U_C(x)\leq C$ for all $x\in\rset$, $U_C(x)=C$ for $x\leq -2$,
$U_C(x)=0$ for $x\geq 0$, and $-C\leq U'_C(x)\leq 0$ for $x\in[-2,0]$.
Note that the potential $U(x,t)$ satisfies \eqref{Bound}.

Let $\gamma^x(t)$, $t\in[t_1,t_2]$, be a minimizer with
\begin{equation}
  \label{defRTrep}
  |\gamma^x(t_2)-y|=|x-y|\leq R_T\equiv K(\log T)^{2/\beta}/2. 
\end{equation}
Without loss of ge\-ne\-rality suppose $T=t_2-t_1>1$ and $y=0$. To establish
Theorem~\ref{lowerbound}, we consider three possible cases:
(i)~$\gamma^x(t_2-1)\leq-g_T(1)$, (ii)~$\gamma^x(t_2-1)>-g_T(1)$ and
$x\geq 0$, and (iii)~$\gamma^x(t_2-1)>-g_T(1)$ and $x<0$. Lemmas
\ref{t:i}--\ref{t:iii} cover each of these cases and together complete the
proof.

\begin{lm}[case (i)]
  \label{t:i}
  If $\gamma^x(t_2-1)\leq-g_T(1)$, then for any $K>0$ there holds
  $\dot\gamma^x(t_2)\geq K(\log T)^{2/\beta}/2^{\beta/(\beta-1)}$ for
  $T$ large enough.
\end{lm}

\begin{proof}
For the average velocity of $\gamma^x$ at the instant~$t_2-1$ we have
\begin{equation}
  \label{temp3}
  w(1)=|x-\gamma^x(t_2-1)|\geq g_T(1)-R_T\geq K(\log T)^{2/\beta}/2,
\end{equation}
where we use inequalities \eqref{residue} and~\eqref{defRTrep}. Thus
the hypothesis of Lemma~\ref{t:w_1} is satisfied if $K(\log
T)^{2/\beta}\geq 2(2C)^{1/(\beta-1)}$, which by the first of
inequalities~\eqref{averagedot} then implies that
$\dot\gamma^x(t_2)\geq2^{-1/(\beta-1)}w(1)$ and, together with
estimate~\eqref{temp3} for~$w(1)$, gives the statement of the
lemma.
\end{proof}

\begin{lm}[case (ii)]
  \label{t:ii}
  Let $\gamma^x(t_2-1)>-g_T(1)$, $x\geq 0$ and
  $K=(C\beta/5)^{1/\beta}$. Then $\dot\gamma^x(t_2)\geq K(\log
  T)^{2/\beta}/2$ for $T$ large enough.
\end{lm}

\begin{proof}
We first note that the minimizer~$\gamma^x$ cannot stay in the domain
where $U=0$ for all $t\in[t_1,t_2]$. More formally, define
\begin{equation}
  \bs\equiv\inf\{\,s\in(1,T)\mid\gamma^x(t_2-s)\leq-g_T(s)\,\};
\end{equation}
then $\bs<T$ and $\gamma^x(t_2-\bs)=-g_T(\bs)$. Indeed, otherwise the
velocity of the minimizer~$\gamma^x$ would vanish for all~$t$ and we
would have $A_{t_1,t_2}[\gamma^x]=0$. Consider a continuous trajectory
$\bar\gamma$ defined on $[t_1,t_2]$ by
\begin{equation}
  \bar\gamma(t_2-s)\equiv
  \begin{cases}
    x-(x+g_T(1)+2)s, & s\in[0,1),\\
    -g_T(s)-2, & s\in[1,T].
  \end{cases}
\end{equation}
Using \eqref{Bound} and~\eqref{integr}, we obtain the following
estimate for the action $A_{t_1,t_2}[\bar\gamma]$:
\begin{equation}
  \label{gamma0}
  \begin{array}{r@{\;}c@{\;}l}
    A_{t_1,t_2}[\bar\gamma]&=&
      A_{t_2-1,t_2}[\bar\gamma]+A_{t_1,t_2-1}[\bar\gamma]\\[1ex]
    &\leq&\displaystyle\frac{(x+g_T(1)+2)^\beta}{\beta}
      +\frac 1\beta\int_0^T(\dot g_T(s))^\beta\,ds-C(T-1)\\[2ex]
    &=&\displaystyle\frac{(x+g_T(1)+2)^\beta}{\beta}+\frac{2K^\beta T}{\beta}-C(T-1).
  \end{array}
\end{equation}
Observing that $A_{t_1,t_2}[\bar\gamma]\geq A_{t_1,t_2}[\gamma^x]=0$
and using the fact that $K^\beta=C\beta/5$, we derive
\begin{equation}
  \frac 35 T\le 1+\frac{(x+g_T(1)+2)^\beta}{\beta C}.
\end{equation}
Since, for $T$ large enough, $x\leq K(\log T)^{2/\beta}/2$ and
$g_T(1)\leq 2K(\log T)^{2/\beta}$ by \eqref{defRTrep}
and~\eqref{residue}, we see that the hypothesis $\bs=T$ leads to a
contradiction.

If $\bs\leq 3$, then the statement of this lemma is established by the
same argument as in Lemma~\ref{t:i}. Therefore assume that
$\gamma^x(t_2-\bs)=-g_T(\bs)$ with $3<\bs<T$ and consider the
continuous trajectory~$\gamma$ defined for $t\in[t_2-\bs,t_2]$ by
\begin{equation}
  \label{defgammaii}
  \gamma(t_2-s)\equiv
  \begin{cases}
    x-(x+g_T(1))s,& s\in[0,1),\\
    -g_T(s)-2(s-1),& s\in[1,2),\\
    -g_T(s)-2,& s\in[2,\bs-1),\\
    -g_T(s)-2(\bs-s),& s\in[\bs-1,\bs].
  \end{cases}
\end{equation}
For the action~$A_{t_2-\bs,t_2}[\gamma]$ we get using~\eqref{Bound} that
\begin{equation}
  \label{actionii}
  A_{t_2-\bs,t_2}[\gamma]
  \leq\frac{(x+g_T(1))^\beta}{\beta}
    +\frac 1\beta\int_1^{\bs}(\dot g_T(s))^\beta\,ds
    +\frac{K^\beta}{\beta}(I_1+I_2)-C(\bs-3),
\end{equation}
where $I_1$ and~$I_2$ are defined by
\begin{equation}
  \label{defI1I2}
  \begin{split}
    I_1\equiv\frac 1{K^\beta}\int_1^2
    \bigl((\dot g_T(s)+2)^\beta-(\dot g_T(s))^\beta\bigr)\,ds,\\
    I_2\equiv\frac 1{K^\beta}\int_{\bs-1}^{\bs}
    \bigl(|\dot g_T(s)-2|^\beta-(\dot g_T(s))^\beta\bigr)\,ds.
  \end{split}
\end{equation}
Note also that by~\eqref{Main}
\begin{equation}
  \label{est4}
  \begin{array}{r@{\;}c@{\;}l}
    \displaystyle\frac1\beta\int_1^{\bs}(\dot g_T(s))^\beta\,ds&<&\displaystyle
      \frac{4K^\beta\bs}{\beta}+\frac{(g_T(\bs))^\beta}{\beta\bs^{\beta-1}}
      -\frac{1}{\beta}\int_0^1(\dot g_T(s))^\beta\,ds\\[2ex]
    &\leq&\displaystyle
      \frac{4K^\beta\bs}{\beta}
      +\frac{(g_T(\bs)-g_T(1))^\beta}{\beta(\bs-1)^{\beta-1}},
  \end{array}
\end{equation}
where the last line follows from Jensen's inequality. On the other
hand, since for $t\in[t_2-\bs,t_2]$ the minimizer~$\gamma^x$ stays in
the domain where $U=0$, is velocity remains constant and we have
\begin{equation}
  \label{temp5}
  A_{t_2-\bs,t_2}[\gamma^x]=\frac{(x+g_T(\bs))^\beta}{\beta\bs^{\beta-1}}
\end{equation}
Plugging \eqref{actionii}, \eqref{est4} and~\eqref{temp5} into the
inequality $A_{t_2-\bs,t_2}[\gamma]-A_{t_2-\bs,t_2}[\gamma^x]\geq 0$
gives
\begin{equation}
  \label{ineqs0}
  \bs<15+I_1+I_2
    +\displaystyle\frac{(x+g_T(1))^\beta}{K^\beta}
    +\frac{(g_T(\bs)-g_T(1))^\beta}{K^\beta(\bs-1)^{\beta-1}}
    -\frac{(x+g_T(\bs))^\beta}{K^\beta\bs^{\beta-1}}.
\end{equation}
where we took into account that $C=5K^\beta\beta^{-1}$.

We now estimate terms in the right-hand side of \eqref{ineqs0}. Note
first that for $T$ large enough
\begin{equation}
  \label{estI1}
  \begin{array}{r@{\;}c@{\;}l}
    I_1&=&\displaystyle\frac 1{K^\beta}\int_1^2
      \left(\bigl(\dot g_T(s)+2\bigr)^\beta
      -(\dot g_T(s))^\beta\right)ds\\
    &<&\displaystyle(\log T)^2\int_1^2
      \left(\left(1-\frac{\log s}{\log T}\right)^{2/\beta}
      +\frac2{K(\log T)^{2/\beta}}\right)ds<2(\log T)^2
  \end{array}
\end{equation}
and similarly $I_2<2(\log T)^2$. Second, note that if $T$ is so large
that the right-hand side of~\eqref{residue} is less than $2K(\log
T)^{2/\beta}$ for $s=1$, then by \eqref{residue} and~\eqref{defRTrep}
\begin{equation}
  \label{temp4}
  \frac{(x+g_T(1))^\beta}{K^\beta}\leq(5/2)^\beta(\log T)^2.
\end{equation}
Third, note that since~$x\geq0$ we can use~\eqref{s_2} to get
\begin{equation}
  \label{temp6}
  \frac{(g_T(\bs)-g_T(1))^\beta}{K^\beta(\bs-1)^{\beta-1}}
    -\frac{(x+g_T(\bs))^\beta}{K^\beta \bs^{\beta-1}}
    <\bar M(\log T)^2.
\end{equation}
Taking the estimates for $I_1$ and~$I_2$ (see~\eqref{estI1}),
\eqref{temp4} and~\eqref{temp6} into account in~\eqref{ineqs0}, we get
$\bs\leq M(\log T)^2$ for $T$ large enough with a suitable
constant~$M$.

Now, using again the fact that the velocity of the
minimizer~$\gamma^x$ stays constant for $t\in[t_2-\bs,t_2]$, we
get from~\eqref{residue} for large enough~$T$ that
\begin{equation}
  \label{finalii}
  \dot\gamma^x(t_2)=\frac{x+g_T(\bs)}{\bs}\geq\frac{g_T(\bs)}{\bs}
  > K(\log T)^{2/\beta}\left(1-\frac{\log \bs}{\log T}\right)^{2/\beta}
  > \frac{K}{2}(\log T)^{2/\beta},
\end{equation}
which establishes the statement of Lemma~\ref{t:ii}.
\end{proof}

\begin{lm}[case (iii)]
  \label{t:iii}
  Let $\gamma^x(t_2-1)>-g_T(1)$, $x<0$ and $K=(C\beta/5)^{1/\beta}$.
  Then $\dot\gamma^x(t_2)\geq K(\log T)^{2/\beta}/2^{\beta/(\beta-1)}$
  for $T$ large enough.
\end{lm}

\begin{proof}
Take $(t_2-\bar s,t_2-\bbs)$ to be the largest neighbourhood of the
instant~$t_2-1$ in which the minimizer~$\gamma^x$ stays in the domain
where $U=0$. More formally, define
\begin{equation}
  \begin{split}
    \bs\equiv\inf\{\,s\in(1,T)\mid\gamma^x(t_2-s)\leq-g_T(s)\,\},\\
    \bbs\equiv\sup\{\,s\in(0,1)\mid\gamma^x(t_2-s)\leq-g_T(s)\,\}.
  \end{split}
\end{equation}

Since $x<0$, the minimizer~$\gamma^x$ must intersect the
curve~$-g_T(t_2-t)$ for $t>t_2-1$, so $\gamma^x(t_2-\bbs)=-g_T(\bbs)$.
Moreover, observe that $\bs<T$ and $\gamma^x(t_2-\bs)=-g_T(\bs)$.
Indeed, otherwise the minimizer~$\gamma^x$ would necessarily stay in
the domain where $U=0$ for all $t\in[t_1,t_2-\bbs]$, so its
velocity would have to vanish and we would have
$\gamma^x(t)=\gamma^x(t_2-\bbs)$ and
$A_{t_1,t_2-\bbs}[\gamma^x]=0$. Assuming, without loss of generality,
that $T\equiv t_2-t_1>\bbs+1$, consider a continuous trajectory~$\bar\gamma$
defined on $[t_1,t_2]$ by
\begin{equation}
  \bar\gamma(t_2-s)\equiv
  \begin{cases}
    -g_T(s)-2(s-\bbs), & s\in[\bbs,\bbs+1),\\
    -g_T(s)-2, & s\in[\bbs+1,T].
  \end{cases}
\end{equation}
Assuming $T$ so large that $\dot g_T(s)=K\bigl(\log(T/s)\bigr)^{2/\beta}>2$ for
$s\in[0,2]$, and using \eqref{Bound}, the inequalities $0\leq\bbs<1$,
and~\eqref{integr}, we obtain the following estimate for the action
$A_{t_1,t_2-\bbs}[\bar\gamma]$:
\begin{equation}
  \begin{array}{r@{\;}c@{\;}l}
    \llap{$A_{t_1,t_2-\bbs}$}[\bar\gamma]&=&\displaystyle
      \int_\bbs^{\bbs+1}\left(
      \frac 1\beta(\dot g_T(s)+2)^\beta-
      U\bigl(\bar\gamma(t_2-s),t_2-s\bigr)\right)ds\\[1ex]
    &&\displaystyle{}+\frac 1\beta\int_{\bbs+1}^{T}(\dot g_T(s))^\beta\,ds
      -C(T-\bbs-1)\\[2ex]
    &\leq&\displaystyle
      \frac 1\beta\int_0^2(\dot g_T(s)+2)^\beta\,ds
      +\frac 1\beta\int_0^{T}(\dot g_T(s))^\beta\,ds-C(T-2)\\[3ex]
    &<&\displaystyle\frac{2(2K)^\beta}{\beta}\bigl((\log\frac T2+1)^2+1\bigr)
      +2C-\left(C-\frac{2K^\beta}{\beta}\right)T.
  \end{array}
\end{equation}
Now note that $\gamma^x$ is a minimizer, so we must have
$A_{t_1,t_2-\bbs}[\bar\gamma]\geq A_{t_1,t_2-\bbs}[\gamma^x]=0$. Thus
the hypothesis $\bs=T$ leads to contradiction, since for
$K=(C\beta/5)^{1/\beta}$ the right-hand side of the last inequality
becomes negative for large~$T$.

If $\bs\leq 3$, then the statement of this lemma is established by the
same argument as in Lemma~\ref{t:i}. Assuming that
$\gamma^x(t_2-\bs)=-g_T(\bs)$ and $\gamma^x(t_2-\bbs)=-g_T(\bbs)$ with
$0\leq \bbs<1$ and $3<\bs<T$, consider the continuous
trajectory~$\gamma$ defined for $t\in[t_2-\bs,t_2-\bbs]$ by
\begin{equation}
  \label{defgammaiii}
  \gamma(t_2-s)\equiv
  \begin{cases}
    -g_T(s)-2(s-\bbs),& s\in[\bbs,\bbs+1),\\
    -g_T(s)-2,& s\in[\bbs+1,\bs-1),\\
    -g_T(s)-2(\bs-s),& s\in[\bs-1,\bs].
  \end{cases}
\end{equation}
Using~\eqref{Bound}, we estimate the action
$A_{t_2-\bs,t_2-\bbs}[\gamma]$ by
\begin{equation}
  \label{actioniii}
    A_{t_2-\bs,t_2}[\gamma]\leq C(\bbs+2-\bs)
      +\frac 1\beta\left(\int_{\bbs}^{\bs}(\dot g_T(s))^\beta\,ds\right)
        +\frac{K^\beta}{\beta}(I_1+I_2),
\end{equation}
where $I_1$ and~$I_2$ are defined by formulas~\eqref{defI1I2}, except
that $I_1$ invovles integration from $\bbs$ to~$\bbs+1$. Note also
that, similarly to~\eqref{est4},
\begin{equation}
  \label{temp7}
  \frac 1\beta\int_\bbs^\bs(\dot g_T(s))^\beta\, ds
  <\frac{4K^\beta\bs}{\beta}
  +\frac{(g_T(\bs)-g_T(\bbs))^\beta}{\beta(\bs-\bbs)^{\beta-1}}
\end{equation}
On the other hand, since for $t\in[t_2-\bs,t_2-\bbs]$ the
minimizer~$\gamma^x$ stays in the domain where $U\equiv 0$, its
velocity remains constant and we have
\begin{equation}
  \label{temp8}
  A_{t_2-\bs,t_2-\bbs}[\gamma^x]=
    \frac{(g_T(\bs)-g_T(\bbs))^\beta}{\beta(\bs-\bbs)^{\beta-1}}.
\end{equation}
Plugging \eqref{actioniii}, \eqref{temp7} and~\eqref{temp8} into the
inequality
$A_{t_2-\bs,t_2-\bbs}[\gamma]-A_{t_2-\bs,t_2-\bbs}[\gamma^x]\geq 0$
and taking into account that $C=5K^\beta\beta^{-1}$, we get a simpler
form of inequality~\eqref{ineqs0}:
\begin{equation}
  \label{ineqs0mod}
  \bs<5(\bbs+2)+I_1+I_2.
\end{equation}
However, this time we need a more accurate estimate of the
sum~$I_1+I_2$ than~\eqref{estI1} can give. Indeed, in the present
case, unlike case~(ii), we have only indirect control over
$\dot\gamma^x(t_2)$, namely that provided by Lemma~\ref{t:w_1}; this
requires a more stringent constraint on~$\bs$.

Recall that \eqref{estI1} tells that $I_1$ and~$I_2$, and therefore $\bs$,
are not larger than $O\bigl((\log T)^2\bigr)$. Thus for suitably large~$T$ we
can expand integrands in $I_1$ and~$I_2$:
\begin{equation}
  \label{I1I2exp}
  \begin{array}{r@{\;}c@{\;}l}
    I_1&=&\displaystyle
      \frac 1{K^\beta}\int_\bbs^{\bbs+1}(\dot g_T(s))^\beta
      \left(\left(1+\frac2{\dot g_T(s)}\right)^\beta-1\right)ds\\[3ex]
    &\leq&\displaystyle
      \int_\bbs^{\bbs+1}\left(
      \frac{2\beta}{K}\bigl(\log(T/s)\bigr)^{2-2/\beta}
      +M_1(K)\bigl(\log(T/s)\bigr)^{2-4/\beta}\right)ds,\\[3ex]
    I_2&\leq&\displaystyle
      \int_{\bs-1}^{\bs}\left(
      -\frac{2\beta}{K}\bigl(\log(T/s)\bigr)^{2-2/\beta}
      +M_1(K)\bigl(\log(T/s)\bigr)^{2-4/\beta}\right)ds,
  \end{array}
\end{equation}
where $M_1(K)$ does not depend on~$T$. 

It is easy to check that for $s$ such that $0<s\leq O\bigl(\log T)^2\bigr)$
\begin{equation}
  \label{proofomitted}
  \int_s^{s+1}\bigl(\log(T/u)\bigr)^{2-2/\beta}\, du=(\log T)^{2-2/\beta}
    +O\bigl(((\log T)^{1-2/\beta}\max\bigl(1,(\log(s+1))^3\bigr)\bigr).
\end{equation}
It follows from~\eqref{proofomitted} that
\begin{equation}
    \int_{\bbs}^{\bbs+1}\bigl(\log(T/s)\bigr)^{2-2/\beta}\, ds-
    \int_{\bs-1}^{\bs}\bigl(\log(T/s)\bigr)^{2-2/\beta}\, ds\leq
    M_2(\log T)^{1-2/\beta}(\log\bs)^3.
\end{equation}

Suppose $1<\beta\leq 2$; then $1-2/\beta\leq 0$ and we have $I_1+I_2\leq
M_2(\log\bs)^3+2M_1$. If $\beta>2$, then $2-4/\beta>1-2/\beta>0$, so that
\begin{equation}
  \label{temp9}
  \int_\bbs^{\bbs+1}\bigl(\log(T/s)\bigr)^{2-4/\beta}\,ds
  <\int_0^1\bigl(\log(T/s)\bigr)^{2-4/\beta}\,ds
  =M_3(\log T)^{2-4/\beta}
\end{equation}
and the rightmost part of~\eqref{temp9} grows with~$T$ faster than
$M_2(\log T)^{1-2/\beta}\log\bs+M_1$; thus
\begin{equation}
  \bs\leq\begin{cases}
    M_2\log\bs+2M_1+5(\bbs+2), &1<\beta\leq 2,\\
    M_3(\log T)^{2-4/\beta}+5(\bbs+2), &\beta>2,
  \end{cases}
\end{equation}
or $\bs\leq M_4(\log T)^{\max\{0,2-4/\beta\}}$ with a suitable
constant~$M_4=M_4(K)$, for large enough~$T$. Note that for such~$\bs$
by \eqref{residue} and~\eqref{defRTrep} we have
\begin{equation}
  w(\bs)=\frac{|x-g_T(\bs)|}{\bs}\geq \frac{K}{2}(\log T)^{2/\beta}.
\end{equation}
Since $2/\beta>\max\{0,2-4/\beta\}/(\beta-1)$ for $\beta>1$, the
condition of Lemma~\ref{t:w_1} is satisfied for large enough~$T$, so
that $|\dot\gamma^x(t_2)|\geq K(\log T)^{2/\beta}/2^{\beta/(\beta-1)}$.
This establishes the statement of Lemma~\ref{t:iii} and concludes the proof of
Theorem~\ref{lowerbound}.
\end{proof}

\subsection{Proof of Theorem~\ref{semiinf}}
In the proof of Theorem~\ref{lowerbound} we constructed an
accelerating potential $U(x,t)$ corresponding to any long enough time
interval $[t_1,t_2]$, $t_2-t_1\equiv T$, and any ball $|x|\leq R_T$ of
terminal positions~$x$ at time~$t_2$. We now glue together a sequence
of such potentials to define for all $t<0$ a potential $U_\infty(x,t)$
that accelerates minimizers indefinitely.

Fix $K=(C\beta/5)^{1/\beta}$. Define increasing sequences $T_n$
and~$S_n$ for $n\geq 1$:
\begin{equation}
  T_1\equiv S_1\equiv\max(1,\bar T),\quad
  T_n\equiv\exp\left(S_{n-1}^{1/\epsilon}\right),\quad
  S_n\equiv S_{n-1}+T_n,\quad n\geq 2,
\end{equation}
where $\bar T$ is large enough so that Theorem~\ref{lowerbound} holds
for $T>\bar T$, and $\epsilon$ is any positive number satisfying
$\epsilon<2(\beta-1)/\beta^2$. Define also
\begin{equation}
  X_0\equiv 0,\quad 
  X_n\equiv\sum_{i=1}^n g_{T_i}(T_i),\quad n\geq 1.
\end{equation}
Note that $g_T(T)=KT\int_0^1|\log x|^{2/\beta}\,dx$ and therefore
$X_n=\bar KS_n$, where $\bar K=K\int_0^1|\log x|^{2/\beta}\,dx$. Finally,
define
\begin{equation}
  U_\infty(x,t)\equiv U_C(x-X_{n-1}+g_{T_n}(-t-S_{n-1}))
\end{equation}
for $t\in(-S_n,-S_{n-1}]$, $n\geq 1$, where $S_0\equiv 0$.

Consider a terminal position~$x$ and take $n$ large enough so that
$|x|\leq\frac12R_{T_n}=\frac K4(\log T_n)^{2/\beta}$. Denote by
$\gamma_n^x(t)$ a minimizer on the time interval $t\in[-S_n,0]$ such
that $\gamma_n^x(0)=x$ and $\dot\gamma_n^x(-S_n)=0$. To
establish Theorem~\ref {semiinf}, we now show that for all~$n$ large
enough
\begin{equation}
  \label{v_n}
  |\dot\gamma_n^x(0)|\geq
  \frac{K_0}{2}(\log T_n)^{\frac2\beta-\epsilon},
\end{equation}
where $K_0=K/2^{2+1/(\beta-1)}$; since $2/\beta-\epsilon>0$, this
implies the statement of the theorem.

To prove \eqref{v_n}, we consider two cases. First assume that
$|\gamma_n^x(-S_{n-1})-X_{n-1}|\leq R_{T_n}=\frac K2(\log
T_n)^{2/\beta}$. Since $\gamma_n^x(t)$ is a minimizer on the time
interval $[-S_n,-S_{n-1}]$ with $\dot\gamma_n^x(-S_n)=0$, it follows
from Theorem~\ref{lowerbound} (with $y=X_{n-1}$) that
$|\dot\gamma_n^x(-S_{n-1})| \geq 2K_0(\log T_n)^{2/\beta}$.
Using~\eqref{EL1} in an argument similar to that of Lemma~\ref{t:w_1},
we obtain
\begin{equation}
  |\dot\gamma_n^x(0)|\geq\left(\left(2K_0(\log T_n)^{2/\beta}\right)^{\beta-1} 
  -CS_{n-1}\right)^{\frac{\scriptstyle1}{\scriptstyle\beta-1}}.
\end{equation}
Observing that $S_{n-1}=(\log T_n)^\epsilon$ and increasing~$n$ if
necessary, we get $|\dot\gamma_n^x(0)|\geq K_0(\log T_n)^{2/\beta}$,
which is even stronger than~\eqref{v_n}.

In the second case, when $|\gamma_n^x(-S_{n-1})-X_{n-1}|>R_{T_n}$,
observe that the average velocity $w(S_{n-1})$ on the interval
$[-S_{n-1},0]$ satisfies the inequality $w(S_{n-1})\geq
\bigl(\frac12R_{T_n}-X_{n-1}\bigr)/S_{n-1}$. Taking into account that
$X_{n-1}=\bar KS_{n-1}$, we obtain for large enough~$n$ that 
\begin{equation}
  w(S_{n-1})\geq \frac{R_{T_n}}{4S_{n-1}}=
  \frac K8(\log T_n)^{\frac2\beta-\epsilon}.
\end{equation}
Using again the facts that $S_{n-1}=(\log T_n)^\epsilon$ and that
$\epsilon< 2(\beta-1)/\beta^2$ and assuming $n$ to be large enough, 
we can ensure that
$w(S_{n-1})>(2CS_{n-1})^{1/(\beta-1)}$. By Lemma~\ref{t:w_1}, this
implies~\eqref{v_n}.
\qed

\section{Conclusion}
\label{s:conclusion}

The results of this paper can be generalized in several directions.
One can consider Lagrangian systems with discrete time. In this
situation one has to find a minimizing sequence
$\{\,x_i\in\rset^d\colon N_1\leq i\leq N_2\,\}$ for the action
\begin{equation}
  A_{N_1,N_2}[\{x_i\}]=
  \sum_{i=N_1}^{N_2-1}\left[\frac{1}{\beta}|x_{i+1}-x_i|^\beta-
  U_i(x_i)\right],
\end{equation}
subject to the condition $x_{N_2}=x$. In physics literature such
systems are called non-stationary Frenkel--Kontorova type models.
Notice that the discrete-time case corresponds to ``kicked forcing''
in the continuous-time setting, i.e., to a forcing of the form
$U(x,t)=\sum_i U_i(x)\delta(t-i)$ (see, e.g., \cite{BFK00}). The results
in the discrete situation are the same as in the continuous-time
setting.

It is also possible to consider more general natural Lagrangian
systems where a Lagrangian has the following form $L(x,v,t)= L_0(v) -
U(x,t)$. This and other generalizations will be discussed in a
forthcoming publication.

It is interesting to study whether in Theorem~\ref{semiinf} it is
possible to replace the one-sided (upper) limit by the two-sided
limit. We believe that the answer to this question is affirmative.

Notice that for the potentials constructed in this paper the partial
derivative $\partial U/\partial t$ is unbounded. It is natural to ask
whether velocity can grow with $T$ in the case when
\begin{equation}
  \left|\frac{\partial U(x,t)}{\partial t}\right|\leq C,\quad 
  x\in\rset^d,t\in\rset.
\end{equation}

It is important to mention that all the ``accelerating'' potentials
constructed in this paper have a very specific form. We expect that
for generic bounded time-dependent potentials the velocity of minimizers
is bounded. Below we formulate this statement as a conjecture in the
case of random potentials.
\begin{conj}
  Let
  \begin{equation}
    U(x,t)=\sum_{j=1}^N U_j(x)a_j^\omega(t),\quad 
    x\in\rset^d,t\in\rset,
  \end{equation} 
  where $U_j(x)$ are fixed non-random potentials of class~$C^1$
  satisfying condition~\eqref{Bound} and $(a_j^\omega(t),1\leq j\leq
  N)$ is a realization of a stationary vector-valued random process
  with exponentially decaying correlation, where
  $\omega$ is a point of the corresponding probability space and
  $\sup_{j,t}|(a_j^\omega(t)|\leq 1$ for almost all~$\omega$. Then
  there exists a random constant $C^\omega(x)$ such that uniformly for all $t\leq -1$
  \begin{equation}
    |\dot \gamma^x_{t,0}(0)| \leq C^\omega(x), 
  \end{equation}
  where $\gamma^x_{t,0}(\tau)$ is a
  minimizer on $[t,0]$ such that $\gamma^x_{t,0}(0)=x$.
\end{conj}

If this conjecture holds true, then global solutions exist with
probability~1 in the case of random potentials.

\subsection*{Acknowldegments}
The work which led to the present paper was mostly carried out during
the stay of A.S.\ in the Isaac Newton Institute in Cambridge on the
INTAS fellowship YSF2001:172 and later when K.Kh.\ and A.S.\ attended
a workshop on Idempotent Mathematics and Mathematical Physics
organized by G.~Litvinov and V.~Maslov at the Erwin Schr{\"o}dinger
Institute in Vienna. We gratefully acknowledge the support and
hospitality of these foundations, institutions, and individuals. A.S.\
also acknowledges the support of the French Ministry of education,
CNRS, and the Russian Foundation for basic research (project
02--01--1062). Finally, it is our pleasant duty to thank A.~Fathi,
S.~Illman, A.~Kelbert, G.~Paternain and A.~Teplinsky for helpful
discussions.

\appendix

\section{Proof of Lemma~\ref{t:g_T}}
\label{a:g_T}

Eq.~\eqref{integr} is obtained by two integrations by parts. 

To obtain~\eqref{residue}, integrate the right-hand side
of~\eqref{defgT} by parts twice to get
\begin{equation}
  \label{temp1}
  \begin{array}{r@{}l}
  g_T(s)={}&\displaystyle Ks\bigl(\log(T/s))^{2/\beta}
    \left(1+\frac2{\beta\log(T/s)}\right)\\
  &\displaystyle{}+\frac{2(2-\beta)K}{\beta^2}
    \int_0^s\bigl(\log(T/u)\bigr)^{(2/\beta)-2}du.
  \end{array}
\end{equation}
Making the change of variable $v=\log(T/u)$ in the integral in the
right-hand side of this formula, we get
\begin{equation}
  \int_0^s\bigl(\log(T/u)\bigr)^{(2/\beta)-2}\,du=
  T\int_{\log(T/s)}^\infty v^{(2/\beta)-2}e^{-v}\,dv\leq
  s\bigl(\log(T/s)\bigr)^{(2/\beta)-2}
\end{equation}
because $v^{(2/\beta)-2}\leq\bigl(\log(T/s)\bigr)^{(2/\beta)-2}$ for
$v\geq\log(T/s)$ and~$\beta>1$. Together with~\eqref{temp1} this
implies~\eqref{residue}. 

We now use~\eqref{residue} to obtain
\begin{equation}
  \begin{array}{r@{\,}c@{\,}l}
    s^{1-\beta}\bigl(g_T(s)\bigr)^\beta&\ge&\displaystyle
    K^\beta s\bigl(\log(T/s)\bigr)^2\left(1+\frac2{\beta\log(T/s)}-
      \frac{2|2-\beta|}{(\beta\log(T/s))^2}\right)^\beta\\[2ex]
    &\ge&\displaystyle
    K^\beta s\bigl(\log(T/s)\bigr)^2\left(1+\frac2{\log(T/s)}-
      \frac{2|2-\beta|}{\beta(\log(T/s))^2}\right);
  \end{array}
\end{equation}
the last line here follows from the inequality $(1+z)^\beta\geq1+\beta
z$ valid for $\beta>1$. Together with~\eqref{integr} this gives
\begin{equation}
  \int_0^s\frac1\beta|\dot g_T(u)|^\beta\,du-
    \frac{s^{1-\beta}}{\beta}\bigl(g_T(s)\bigr)^\beta\leq
  \frac{2K^\beta s}{\beta}\left(1+\frac{|2-\beta|}{\beta}\right),
\end{equation}
which implies inequality~\eqref{Main} for $\beta>1$.

We finally notice that monotonicity of $g_T(\cdot)$ implies that for
$s>1$
\begin{equation}
  \label{temp2}
  \displaystyle
  \frac{(g_T(s)-g_T(1))^\beta}{(s-1)^{\beta-1}}-
    \frac{(g_T(s))^\beta}{s^{\beta-1}}\leq 
    \left(\frac{g_T(s)}{s}\right)^\beta 
      s\left(\frac1{(1-s^{-1})^{\beta-1}}-1\right).
\end{equation}
Observe that for $x\in[0,1/3]$
\begin{equation}
  \label{temp10}
  \frac1{(1-x)^{\beta-1}}-1\leq 
  3\left(\biggl(\frac32\biggr)^{\beta-1}-1\right)x
\end{equation}
(note that the left-hand side of \eqref{temp10} is a
convex function, whose graph on the specified interval lies below its chord given by
the right-hand side). Furthermore, notice that
\begin{equation}
  \label{temp12}
  \frac{g_T(s)}{s}=\frac Ks\int_0^s\bigl(\log(T/u)\bigr)^{2/\beta}\,du
  =K\int_0^\infty\bigl(v+\log(T/s)\bigr)^{2/\beta}\mathrm{e}^{-v}\,dv,
\end{equation}
where we performed the change of variable $v=\log(s/u)$. If $T/2<s\leq
T$, then the right-hand side of this expression is bounded uniformly
in~$T$; for $1<s\leq T/2$ we have
\begin{equation}
  \begin{array}{l}
    \displaystyle K\int_0^\infty\bigl(v+\log(T/s)\bigr)^{2/\beta}\mathrm{e}^{-v}\,dv=
    K\bigl(\log(T/s)\bigr)^{2/\beta}\int_0^\infty
      \left(\frac v{\log(T/s)}+1\right)^{2/\beta}\mathrm{e}^{-v}\,dv\\
    \displaystyle\quad\leq K\bigl(\log(T/s)\bigr)^{2/\beta}\int_0^\infty
      \left(\frac v{\log 2}+1\right)^{2/\beta}\mathrm{e}^{-v}\,dv.
  \end{array}
\end{equation}
Therefore for $T$ large enough
\begin{equation}
  \label{temp11}
  \frac{g_T(s)}{s}\leq\tilde K\bigl(\log(T/s)\bigr)^{2/\beta}
\end{equation}
with a suitable $\tilde K>0$. Inequalities \eqref{temp2},
\eqref{temp10}, and~\eqref{temp11} together give~\eqref{s_2} for
$3\leq s\leq T$.

\section{Absense of blow-up in the time-periodic case}
\label{a:fathi}

In this appendix we present A.~Fathi's proof that there is no blow-up
if the potential $U(x,t)$ is periodic in time. Therefore, in addition
to assumptions~\eqref{Bound}, we require that $U(x,t)=U(x,t+1)$ for any
$x\in\rset^d$ and any $t\in\rset$.

Let $x,y\in\rset^d$, $t_1<t_2$. Since the action
functional~\eqref{defaction} is bounded below, we can write, repeating
definition~\eqref{defkernel},
\begin{equation}
  A_{t_1,t_2}(y,x)=\inf_{\gamma(t_1)=y,\;\gamma(t_2)=x}A_{t_1,t_2}[\gamma].
\end{equation}
In what follows we assume that this infimum is attained, which is a standard
result under the present hypotheses on the Lagrangian (see, e.g.,
\cite{F01}). The following elementary lemma is also standard.

\begin{lm}
  \label{t:Alip}
  The function $A_{t_1,t_2}(y,x)$ is uniformly locally Lipschitz: for
  any $W>0$, there exists $K=K(W,t_1,t_2)$ such that if $t_1<t_2$ and
  $x_1,x_2,y\in\rset^d$ are such that $|x_i-y|\leq W\cdot(t_2-t_1)$,
  $i=1,2$, then
  \begin{equation}
    \label{alip0}
    |A_{t_1,t_2}(y,x_2)-A_{t_1,t_2}(y,x_1)|\leq K|x_2-x_1|
  \end{equation}
  Moreover, the function $A_{t_1,t_2}(y,x)$ admits the following
  bounds: for any $x,y\in\rset^d$, $t_1<t_2$
  \begin{equation}
  \label{alip1}
    \frac{|x-y|^\beta}{\beta(t_2-t_1)^\beta}-C\leq
    \frac1{t_2-t_1}A_{t_1,t_2}(y,x)
    \leq\frac{|x-y|^\beta}{\beta(t_2-t_1)^\beta}.
  \end{equation}
\end{lm}

\begin{proof}
Let $\gamma_0$ be a minimizing curve 
%for which $A_{t_1,t_2}[\gamma_0]=A_{t_1,t_2}(\gamma_0(t_1),\gamma_0(t_2))$ 
and $w={|\gamma_0(t_2)-\gamma_0(t_1)|}/(t_2-t_1)$ be its average
velocity defined as in~\eqref{defws} above. By classic arguments, the
Lipschitz property of $A_{t_1,t_2}(y,x)$ follows from boundedness of
$|\dot\gamma(t)|$ on $[t_1, t_2]$, which itself is established in a way
similar to Lemma~\ref{t:w_1}.

The left inequality in \eqref{alip1} follows from Jensen's inequality
\eqref{estimate} and condition~\eqref{Bound}. The right inequality
follows in a similar way from the inequality
$A_{t_1,t_2}[\gamma_0]\leq A_{t_1,t_2}[\gamma]$ written for
$\gamma(t)=\frac{t_2-t}{t_2-t_1}y+ \frac{t-t_1}{t_2-t_1}x$.
\end{proof}

Following A.~Fathi, we introduce two concepts now. A function
$S\colon\rset^d\to\rset$ is said to be $(L,t_1,t_2)$-dominated for a
time interval $[t_1,t_2]$ and a constant $L\in\rset$, if for any
$x,y\in\rset^d$
\begin{equation}
  \label{defdomin}
  S(x)-S(y)\leq A_{t_1,t_2}(y,x)+L(t_2-t_1),
\end{equation}
and Lipschitz in the large with constant~$K$, if for any $x,y\in\rset^d$ 
\begin{equation}
  \label{defliplarge}
  |S(x)-S(y)|\leq K(|x-y|+1)
\end{equation}
with some~$K>0$.

\begin{lm}
  \label{t:Sliplarge}
  An $(L,t_1,t_2)$-dominated function is Lipschitz in the large with
  constant~$K$ depending on $L$ and~$t_1,t_2$.
\end{lm}

\begin{proof}
For $x,y\in\rset^d$ define the sequence $(x_i)$, $0\leq i\leq[|x-y|]$,
by $x_i=y+ir$, where $r=|x-y|^{-1}(x-y)$ is a unit vector collinear
with $x-y$. (Here, as above, $[\cdot]$ stands for the integer part.)
We can write
\begin{equation}
  S(x)-S(y)=\sum_{1\leq i\leq[|x-y|]}(S(x_i)-S(x_{i-1}))
    +S(x)-S(x_{[|x-y|]}).
\end{equation}
Using the property of $(L,t_1,t_2)$-domination and the right
inequality~\eqref{alip1}, we get
\begin{equation}
  S(x)-S(y)\leq
    \left(\frac1{\beta(t_2-t_1)^{\beta-1}}+L(t_2-t_1)\right)
    \bigl([|x-y|]+1\bigr).
\end{equation}
Together with the reverse inequality obtained by interchanging the
roles of $x$ and~$y$, this implies~\eqref{defliplarge}.
\end{proof}

Denote the Lax--Ole{\u\i}nik solution operator over a time interval
$[t_1,t_2]$ for the Cauchy problem for equation \eqref{HJE} by
\begin{equation}
  \label{defT}
  T_{t_1,t_2}S(x)\equiv\inf_{y\in\rset^d}(A_{t_1,t_2}(y,x)+S(y)).
\end{equation}

\begin{lm}
  \label{t:Tdomin}
  For any $L$ and any $t_1<t_2$, the operator $T_{t_1,t_2}$ maps the
  set of $(L,t_1,t_2)$-dominated functions into itself.
\end{lm}

\begin{proof}
If $S(x)$ is $(L,t_1,t_2)$-dominated, it follows from~\eqref{defdomin} that for
any $x\in\rset^d$
\begin{equation}
  S(x)\leq\inf_{y\in\rset^d}(A_{t_1,t_2}(y,x)+S(y)+L(t_2-t_1))=
  T_{t_1,t_2}S(x)+L(t_2-t_1).
\end{equation}
Therefore for any $z\in\rset^d$
\begin{equation}
  \begin{split}
    T_{t_1,t_2}S(x)&=\inf_{y\in\rset^d}(A_{t_1,t_2}(y,x)+S(y))\\
    &\leq A_{t_1,t_2}(z,x)+S(z)\leq A_{t_1,t_2}(z,x)+T_{t_1,t_2}S(z)+L(t_2-t_1),
  \end{split}
\end{equation}
which implies $(L,t_1,t_2)$-domination for $T_{t_1,t_2}S(x)$.
\end{proof}

\begin{lm}
  \label{t:Tlip}
  For any $K>0$ and any $t_1<t_2$, the operator $T_{t_1,t_2}$ maps the
  set of functions that are Lipschitz in the large with constant $K$
  into the set of Lipschitz functions with constant $\bar K=\bar
  K(K,t_1,t_2)$.
\end{lm}

\begin{proof}
Let $S(x)$ be a function that is Lipschitz in the large with
constant~$K$. Then for any $y\in\rset^d$
\begin{equation}
    A_{t_1,t_2}(y,x)+S(y)\geq A_{t_1,t_2}(y,x)+S(x)-K(|x-y|+1).
\end{equation}
On the other hand, by definition~\eqref{defT} of the
operator~$T_{t_1,t_2}$ and the last inequality in~\eqref{alip1}, we
have
\begin{equation}
  T_{t_1,t_2}S(x)\leq A_{t_1,t_2}(x,x)+S(x)\leq S(x).
\end{equation}
Therefore instead of~\eqref{defT} we can write
\begin{equation}
  \label{ersatzT}
  T_{t_1,t_2}S(x)=\inf_{y\in\rset^d}
    (A^K_{t_1,t_2}(y,x)+S(y)),
\end{equation}
where $A^K_{t_1,t_2}(y,x)=\min\{A_{t_1,t_2}(y,x),K(|x-y|+1)\}$. 

The first inequality in~\eqref{alip1} implies that
$A^K_{t_1,t_2}(y,x)=K(|x-y|+1)$ if $|x-y|>R$ with a suitable
$R=R(K,t_1,t_2)$. Together with the first part of Lemma~\ref{t:Alip}
this means that $A^K_{t_1,t_2}(y,x)$ is a Lipschitz function of~$x$,
with a constant~$\bar K=\bar K(K,t_1,t_2)$ that does not depend
on~$y$. It now follows from~\eqref{ersatzT} that $T_{t_1,t_2}S(x)$ is
Lipschitz with the same constant.
\end{proof}

Now observe that by \eqref{alip1} any constant function is
$(L,t_1,t_2)$-dominated with $L=C$ for any $t_1<t_2$. Using periodicity of $U$ and Lemmas~\ref{t:Sliplarge}--\ref{t:Tlip}
with $t_1=n$, $t_2=n+1$ for integer $n\geq 0$, we see that the
solution~$S(x,t)$ of the Cauchy problem for equation~\eqref{HJE} with
the initial condition $S(x,0)=0$ stays $(C,0,1)$-dominated and therefore
Lipschitz for all integer moments of time. Applying, for any
noninteger $t>0$, Lemma~\ref{t:Tlip} again with $t_1=[t]$, $t_2=t$, we
get Lipschitzness for all $t>0$ with a suitable constant depending on
the parameters of the problem.

%\bibliographystyle{amsalpha}
%\bibliography{velocity}

\begin{thebibliography}{EKMS00}

\bibitem[BFK00]{BFK00}
J.~Bec, U.~Frisch, and K.~Khanin, \emph{Kicked {B}urgers turbulence}, J. Fluid
  Mech. \textbf{416} (2000), 239--267. \MR{2001d:76071}

\bibitem[CL83]{CL83}
Michael~G. Crandall and Pierre-Louis Lions, \emph{Viscosity solutions of
  {H}amilton-{J}acobi equations}, Trans. Amer. Math. Soc. \textbf{277} (1983),
  no.~1, 1--42. \MR{85g:35029}

\bibitem[EG02]{EG02}
L.~C. Evans and D.~Gomes, \emph{Linear programming interpretations of
  {M}ather's variational principle}, ESAIM Control Optim. Calc. Var. \textbf{8}
  (2002), 693--702 (electronic). \MR{2003h:90032}

\bibitem[EKMS00]{EKMS00}
Weinan E, K.~Khanin, A.~Mazel, and Ya. Sinai, \emph{Invariant measures for
  {B}urgers equation with stochastic forcing}, Ann. of Math. (2) \textbf{151}
  (2000), no.~3, 877--960. \MR{2002e:37134}

\bibitem[Fat01]{F01}
Albert Fathi, \emph{Weak {KAM} theorem in {L}agrangian dynamics}, to be
  published by the Cambridge University Press, 2001.

\bibitem[GIKP03]{GIKP03}
D.~Gomes, R.~Iturriaga, K.~Khanin, and P.~Padilla, \emph{Viscosity limit of
  stationary distributions for the random forced {B}urgers equation}, to be
  published, 2003.

\bibitem[HK03]{HK03}
Viet~Ha Hoang and Konstantin Khanin, \emph{Random {B}urgers equation and
  {L}agrangian systems in non-compact domains}, Nonlinearity \textbf{16}
  (2003), no.~3, 819--842. \MR{1 975 784}

\bibitem[Hop50]{H50}
Eberhard Hopf, \emph{The partial differential equation {$u\sb t+uu\sb x=\mu
  u\sb {xx}$}}, Comm. Pure Appl. Math. \textbf{3} (1950), 201--230.
  \MR{13,846c}

\bibitem[IK03]{IK03}
R.~Iturriaga and K.~Khanin, \emph{Burgers turbulence and random {L}agrangian
  systems}, Comm. Math. Phys. \textbf{232} (2003), no.~3, 377--428. \MR{1 952
  472}

\bibitem[KM97]{KM97}
Vassili~N. Kolokoltsov and Victor~P. Maslov, \emph{Idempotent analysis and its
  applications}, Mathematics and its Applications, vol. 401, Kluwer Academic
  Publishers Group, Dordrecht, 1997. \MR{1 447 629}

\bibitem[Lio82]{L82}
Pierre-Louis Lions, \emph{Generalized solutions of {H}amilton-{J}acobi
  equations}, Research Notes in Mathematics, vol.~69, Pitman (Advanced
  Publishing Program), Boston, Mass., 1982. \MR{84a:49038}

\bibitem[Mat89]{M89}
John~N. Mather, \emph{Minimal measures}, Comment. Math. Helv. \textbf{64}
  (1989), no.~3, 375--394. \MR{90f:58067}

\bibitem[Rou]{R03}
I.V. Roublev, \emph{On two notions of generalized solution to the
  {H}amilton--{J}acobi equation}, Idempotent mathematics and mathematical
  physics.

\bibitem[Sub95]{S95}
Andre\u{\i}~I. Subbotin, \emph{Generalized solutions of first-order {PDE}s},
  Systems \& Control: Foundations \& Applications, Birkh\"auser Boston Inc.,
  Boston, MA, 1995. \MR{96b:49002}

\end{thebibliography}

\providecommand{\bysame}{\leavevmode\hbox to3em{\hrulefill}\thinspace}
\providecommand{\MR}{\relax\ifhmode\unskip\space\fi MR }
\providecommand{\MRhref}[2]{%
  \href{http://www.ams.org/mathscinet-getitem?mr=#1}{#2}
}
\providecommand{\href}[2]{#2}

\end{document}